\documentclass{amsart}%
\usepackage{amssymb}
\usepackage{amsmath}
\usepackage{graphicx}
\usepackage{amsfonts}%
\newtheorem{theorem}{Theorem}
\theoremstyle{plain}

\newtheorem{definition}{Definition}
\newtheorem{example}{Example}

\numberwithin{equation}{section}

\begin{document}
\title{Structure Relations in Special $A_{\infty}$-bialgebras }
\author{Ronald Umble$^{1}$}
\address{Department of Mathematics\\
Millersville University of Pennsylvania\\
Millersville, PA. 17551}
\thanks{$^{1}$ This research was funded in part by a Millersville University faculty
research grant.}
\date{June 22, 2005}

\begin{abstract}
We compute the structure relations in special $A_{\infty}$-bialgebras whose
operations are limited to those defining the underlying $A_{\infty}%
$-(co)algebra substructure. Such bialgebras appear as the homology of certain
loop spaces. Whereas structure relations in general $A_{\infty}$-bialgebras
depend upon the combinatorics of permutahedra, only Stasheff's associahedra
are required here.

\end{abstract}
\maketitle

\section{Introduction}

A general $A_{\infty}$-infinity bialgebra is a DG module $\left(  H,d\right)
$ equipped with a family of structurally compatible operations $\omega
_{j,i}:H^{\otimes i}\rightarrow H^{\otimes j},\,$where $i,j\geq1$ and
$i+j\geq3$ (see \cite{SU3}). In \emph{special} $A_{\infty}$-bialgebras,
$\omega_{j,i}=0$ whenever $i,j\geq2,$ and the remaining operations
$m_{i}=\omega_{1,i}$ and $\Delta_{j}=\omega_{j,1}$ define the underlying
$A_{\infty}$-(co)algebra substructure. Thus special $A_{\infty}$-bialgebras
have the form $\left(  H,d,m_{i},\Delta_{j}\right)  _{i,j\geq2}$ subject to
the appropriate structure relations involving $d,$ the $m_{i}$'s and
$\Delta_{j}$'s. These relations are much easier to describe than those in the
general case, which require the S-U diagonal $\Delta_{P}$ on permutahedra.
Instead, the S-U diagonal $\Delta_{K}$ on Stasheff's associahedra $K=\sqcup
K_{n}$ is required here (see \cite{SU2}).

$A_{\infty}$-bialgebras are fundamentally important structures in algebra and
topology. In general, the homology of every $A_{\infty}$-bialgebra inherits an
$A_{\infty}$-bialgebra structure \cite{SU4}; in particular, this holds for the
integral homology of a loop space. In fact, over a field, the $A_{\infty}%
$-bialgebra structure on the homology of a loop space specializes to the
$A_{\infty}$-(co)algebra structures observed by Gugenheim \cite{Gugenheim} and
Kadeishvili \cite{Kadeishvili1}.

The main result of this paper is the following simple formulation of the
structure relations in special $A_{\infty}$-bialgebras that do not involve
$d$: Let $TH$ denote the tensor module of $H$ and let $e^{n-2}$ denote the top
dimensional face of $K_{n}.$ There is a \textquotedblleft fraction
product\textquotedblright\ on $M=End\left(  TH\right)  $ (denoted here by
\textquotedblleft$\bullet$\textquotedblright) and certain cellular cochains
$\xi,\zeta\in C^{\ast}\left(  K;M\right)  $ such that for each $i,j\geq2,$
\[
\Delta_{j}\bullet m_{i}=\xi^{j}\left(  e^{i-2}\right)  \bullet\zeta^{i}\left(
e^{j-2}\right)  ,
\]
where the exponents indicate certain $\Delta_{K}$-cup powers.

I must acknowledge the fact that many of the ideas in this paper germinated
during conservations with Samson Saneblidze, whose openness and encouragement
led to this paper. For this I express sincere thanks.

\section{Matrix Considerations}

We begin with a brief review of the algebraic machinery we need; for a
detailed exposition see \cite{SU3}. Let $M=\bigoplus\nolimits_{m,n\geq
1}M_{n,m}$ be a bigraded module over a commutative ring $R$ with identity
$1_{R}$ and consider the module $TTM$ of tensors on $TM.$ Given matrices
$X=\left[  x_{ij}\right]  $ and $Y=\left[  y_{ij}\right]  \in\mathbb{N}%
^{q\times p},$ $p,q\geq1,$ consider the submodule%
\begin{align*}
{M}_{Y,X}  &  =\left(  M_{y_{11},x_{11}}\otimes\cdots\otimes M_{y_{1p},x_{1p}%
}\right)  \otimes\cdots\otimes\left(  M_{y_{q1},x_{q1}}\otimes\cdots\otimes
M_{y_{qp},x_{qp}}\right) \\
&  \subset\left(  M^{\otimes p}\right)  ^{\otimes q}\subset TTM.
\end{align*}
Represent a monomial $A=\left(  \theta_{y_{11},x_{11}}\otimes\cdots
\otimes\theta_{y_{1p},x_{1p}}\right)  \otimes\cdots\otimes\left(
\theta_{y_{q1},x_{q1}}\otimes\cdots\otimes\theta_{y_{qp},x_{qp}}\right)
$\linebreak$\in{M}_{Y,X}$ as the $q\times p$ matrix $\left[  A\right]
=\left[  a_{ij}\right]  $ with $a_{ij}=\theta_{y_{ij},x_{ij}}.$ Then $A$ is
the $q$-fold tensor product of the rows of $[A]$ thought of as elements of
$M^{\otimes p};$ we refer to $A$ as a $q\times p$ monomial and often write $A$
when we mean $\left[  A\right]  .$ The \emph{matrix submodule of }$TTM$ is the
sum
\[
\overline{\mathbf{M}}=\bigoplus_{\substack{X,Y\in\mathbb{N}^{q\times
p}\\p,q\geq1}}M_{Y,X}=\bigoplus_{p,q\geq1}(M^{\otimes p})^{\otimes q}.
\]
Given $\mathbf{x\times y=}\left(  x_{1},\ldots,x_{p}\right)  \times\left(
y_{1},\ldots,y_{q}\right)  \in\mathbb{N}^{p}\times\mathbb{N}^{q},$ \ set
$X=\left[  x_{ij}=x_{j}\right]  _{1\leq i\leq q},$ $Y=\left[  y_{ij}%
=y_{i}\right]  _{1\leq j\leq p}$ and denote $\mathbf{M}_{\mathbf{x}%
}^{\mathbf{y}}=M_{Y,X}.$ The \emph{essential submodule of }$TTM$ is
\[
\mathbf{M}=\bigoplus_{\substack{\mathbf{x\times y}\in\mathbb{N}^{p}%
\times\mathbb{N}^{q}\\p,q\geq1}}\mathbf{M}_{\mathbf{x}}^{\mathbf{y}}%
\]
and a $q\times p$ monomial $A\in\mathbf{M}$ has the form%
\[
A=\left[
\begin{array}
[c]{lll}%
\theta_{y_{1},x_{1}} & \cdots & \theta_{y_{1},x_{p}}\\
\multicolumn{1}{c}{\vdots} & \multicolumn{1}{c}{} & \multicolumn{1}{c}{\vdots
}\\
\theta_{y_{q},x_{1}} & \cdots & \theta_{y_{q},x_{p}}%
\end{array}
\right]  .
\]
Graphically represent $A=\left[  \theta_{y_{j},x_{i}}\right]  \in
\mathbf{M}_{\mathbf{x}}^{\mathbf{y}}$ two ways: (1) as a matrix of
\textquotedblleft double corollas\textquotedblright\ in which $\theta
_{y_{j},x_{i}}$ is pictured as two corollas joined at the root--one opening
downward with $x_{i}$ leaves and the other opening upward with $y_{j}$
leaves--and (2) as an arrow in the positive integer lattice $\mathbb{N}^{2}$
from $\left(  \left\vert \mathbf{x}\right\vert ,q\right)  $ to $\left(
p,\left\vert \mathbf{y}\right\vert \right)  ,$ where $\left\vert
\mathbf{u}\right\vert =u_{1}+\cdots+u_{k}$ (see Figure 1).\vspace*{0.4in}

\noindent\hspace*{1.8in}\setlength{\unitlength}{0.001in}\begin{picture}
(152,-130)(202,-90) \thicklines
\put(50,-150){\line( -1,1){150}}
\put(50,-150){\line( 1,1){150}}
\put(50,-150){\line( 0,-1){265}}
\put(50,-150){\makebox(0,0){$\bullet$}}
\put(50,-750){\line( -1,1){150}}
\put(50,-750){\line( 1,1){150}}
\put(50,-750){\line( 1,3){50}}
\put(50,-750){\line( -1,3){50}}
\put(50,-750){\line( 0,-1){265}}
\put(50,-750){\makebox(0,0){$\bullet$}}
\put(550,-150){\line( -1,1){150}}
\put(550,-150){\line( 1,1){150}}
\put(550,-250){\line( -1,-1){150}}
\put(550,-250){\line( 1,-1){150}}
\put(550,-400){\line( 0,1){265}}
\put(550,-150){\makebox(0,0){$\bullet$}}
\put(550,-250){\makebox(0,0){$\bullet$}}
\put(550,-750){\line( -1,1){150}}
\put(550,-750){\line( 1,1){150}}
\put(550,-850){\line( -1,-1){150}}
\put(550,-850){\line( 1,-1){150}}
\put(550,-750){\line( 0,-1){265}}
\put(550,-750){\line( 1,3){50}}
\put(550,-750){\line( -1,3){50}}
\put(550,-750){\makebox(0,0){$\bullet$}}
\put(550,-850){\makebox(0,0){$\bullet$}}
\put(-300,150){\line( 0,-1){1300}}
\put(-300,150){\line( 1,0){100}}
\put(-300,-1150){\line( 1,0){100}}
\put(900,150){\line( 0,-1){1300}}
\put(900,150){\line( -1,0){100}}
\put(900,-1150){\line( -1,0){100}}
\put(-900,-500){\makebox(0,0){$A \in \mathbf M_{1,3}^{2,4} \hspace*{0.1in}
\leftrightarrow  \hspace*{0.1in} $}}
\put(1200,-500){\makebox(0,0){$  \leftrightarrow  \hspace*{0.1in} $}}
\put(1600,-1150){\line(1,0){1100}}
\put(1600,-1150){\line(0,1){1300}}
\put(1800,-1183){\line(0,1){75}}
\put(2200,-1183){\line(0,1){75}}
\put(1563,-950){\line(1,0){75}}
\put(1563,-150){\line(1,0){75}}
\put(1800,-150){\makebox(0,0){$\bullet$}}
\put(2200,-950){\makebox(0,0){$\bullet$}}
\put(2150,-830){\vector(-1,2){300}}
\put(1450,-950){\makebox(0,0){$2$}}
\put(1450,-150){\makebox(0,0){$6$}}
\put(1800,-1300){\makebox(0,0){$2$}}
\put(2200,-1300){\makebox(0,0){$4$}}
\end{picture}\vspace*{1.3in}

\begin{center}
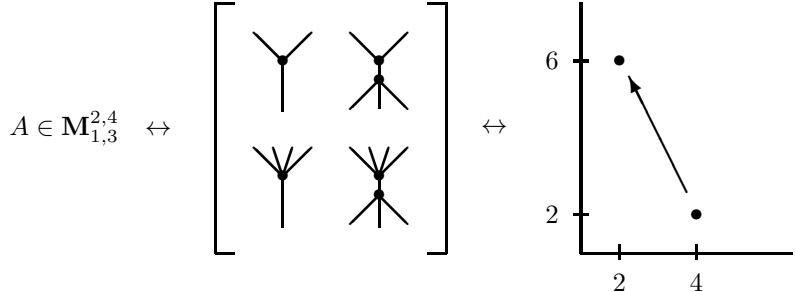
Figure 1. Graphical representations of a typical monomial.\vspace*{0.1in}
\end{center}

Each pairing $\gamma:\bigoplus\nolimits_{r,s\geq1}M^{\otimes r}\otimes
M^{\otimes s}\rightarrow M$ induces an \emph{upsilon product }$\Upsilon
:\overline{\mathbf{M}}\mathbf{\otimes}\overline{\mathbf{M}}\mathbf{\rightarrow
}\overline{\mathbf{M}}$ supported on \textquotedblleft block transverse
pairs,\textquotedblright\ which we now describe.

\begin{definition}
A monomial pair $A^{q\times s}\otimes B^{t\times p}=\left[  \theta_{y_{k\ell
},v_{k\ell}}\right]  \otimes\left[  \eta_{u_{ij},x_{ij}}\right]  \in
\overline{\mathbf{M}}\mathbf{\otimes}\overline{\mathbf{M}}$ is a

\begin{enumerate}
\item[\textit{(i)}] \underline{Transverse Pair} (TP) if $s=t=1,$ $u_{1,j}=q$
and $v_{k,1}=p$ for all $j,k,$ i.e., setting $x_{j}=x_{1,j}$ and
$y_{k}=y_{k,1}$ gives
\[
A\otimes B=\left[
\begin{array}
[c]{c}%
\theta_{y_{1},p}\\
\vdots\\
\theta_{y_{q},p}%
\end{array}
\right]  \otimes\left[
\begin{array}
[c]{lll}%
\eta_{q,x_{1}} & \cdots & \eta_{q,x_{p}}%
\end{array}
\right]  \in\mathbf{M}_{p}^{\mathbf{y}}\otimes\mathbf{M}_{\mathbf{x}}^{q}.
\]

\item[\textit{(ii)}] \underline{Block Transverse Pair} (BTP) if there exist
$t\times s$ block decompositions $A=\left[  A_{k^{\prime}\ell}^{\prime
}\right]  $ and $B=\left[  B_{ij^{\prime}}^{\prime}\right]  $ such that
$A_{i\ell}^{\prime}\otimes B_{i\ell}^{\prime}$ is a TP for\textit{\ all
}$i,\ell$.
\end{enumerate}
\end{definition}

\noindent Unlike the blocks in a standard block matrix, the blocks $A_{i\ell
}^{\prime}$ (or $B_{i\ell}^{\prime}$) in a general BTP may vary in length
within a given row (or column). However, when $A\otimes B\in\mathbf{M}%
_{\mathbf{v}}^{\mathbf{y}}\otimes\mathbf{M}_{\mathbf{x}}^{\mathbf{u}}$ is a
BTP with $\mathbf{u=}\left(  q_{1},\ldots q_{t}\right)  ,$ $\mathbf{v}=\left(
p_{1},\ldots,p_{s}\right)  ,$ $\mathbf{x}=\left(  \mathbf{x}_{1}%
,\ldots,\mathbf{x}_{s}\right)  $ and $\mathbf{y}=\left(  \mathbf{y}_{1}%
,\ldots,\mathbf{y}_{t}\right)  $, the TP $A_{i\ell}^{\prime}\otimes B_{i\ell
}^{\prime}\in\mathbf{M}_{p_{\ell}}^{\mathbf{y}_{i}}\otimes\mathbf{M}%
_{\mathbf{x}_{\ell}}^{q_{i}}$ so that for a fixed $i\,$(or $\ell$) the blocks
$A_{i\ell}^{\prime}$ (or $B_{i\ell}^{\prime}$) have constant length $q_{i}$
(or $p_{\ell}$); furthermore, $A\otimes B$ is a BTP if and only if
$\mathbf{y}\in\mathbb{N}^{|\mathbf{u}|}$ and $\mathbf{x}\in\mathbb{N}%
^{|\mathbf{v}|}$ if and only if the initial point of arrow $A$ coincides with
the terminal point of arrow $B$. Note that BTP block decomposition is unique.

\begin{example}
\label{example1}A pairing of monomials $A^{4\times2}\otimes B^{2\times3}%
\in\mathbf{M}_{2,1}^{1,5,4,3}\otimes\mathbf{M}_{1,2,3}^{3,1}$ is a $2\times2$
BTP per the block decompositions\vspace{0.1in}\newline\hspace*{0.9in}%
\setlength{\unitlength}{.06in}\linethickness{0.4pt}\begin{picture}(118.66,29.34)
\put(5.23,25.00){\makebox(0,0)[cc]{$\theta_{1,2}$}}
\put(5.23,18.00){\makebox(0,0)[cc]{$\theta_{5,2}$}}
\put(5.23,11.67){\makebox(0,0)[cc]{$\theta_{4,2}$}}
\put(5.23,4.67){\makebox(0,0)[cc]{$\theta_{3,2}$}}
\put(1.66,8.67){\dashbox{0.67}(6.33,20.00)[cc]{}}
\put(1.66,2.34){\dashbox{0.67}(6.33,4.67)[cc]{}}
\put(14.00,25.00){\makebox(0,0)[cc]{$\theta_{1,1}$}}
\put(14.00,18.00){\makebox(0,0)[cc]{$\theta_{5,1}$}}
\put(14.00,11.67){\makebox(0,0)[cc]{$\theta_{4,1}$}}
\put(14.00,4.67){\makebox(0,0)[cc]{$\theta_{3,1}$}}
\put(10.33,8.67){\dashbox{0.67}(6.33,20.00)[cc]{}}
\put(10.33,2.34){\dashbox{0.67}(6.33,4.67)[cc]{}}
\put(1.67,15.50){\oval(3.33,28.00)[l]}
\put(16.80,15.50){\oval(3.33,28.00)[r]} \
\put(22.83,16.77){\makebox(0,0)[cc]{\textit{and}}} \
\put(33.00,20.01){\makebox(0,0)[cc]{$\eta_{3,1}$}}
\put(39.58,20.01){\makebox(0,0)[cc]{$\eta_{3,2}$}}
\put(47.53,20.01){\makebox(0,0)[cc]{$\eta_{3,3}$}}
\put(29.67,17.67){\dashbox{0.67}(12.67,5.67)[cc]{}}
\put(44.33,17.67){\dashbox{0.67}(6.67,5.67)[cc]{}}
\put(33.00,12.34){\makebox(0,0)[cc]{$\eta_{1,1}$}}
\put(39.58,12.34){\makebox(0,0)[cc]{$\eta_{1,2}$}}
\put(47.53,13.34){\makebox(0,0)[cc]{$\eta_{1,3}$}}
\put(29.67,10.01){\dashbox{0.67}(12.67,5.67)[cc]{}}
\put(44.53,10.01){\dashbox{0.67}(6.67,5.67)[cc]{}}
\put(28.90,16.27){\oval(2.00,15.33)[l]}
\put(51.56,16.27){\oval(2.00,15.33)[r]}
\put(56.00,16.00){\makebox(0,0)[cc]{$.$}}
\end{picture}\vspace{0.1in}
\end{example}

Given a pairing $\gamma=%
{\textstyle\sum\nolimits_{\mathbf{x\times y}}}
\gamma_{\mathbf{x}}^{\mathbf{y}}:\mathbf{M}_{p}^{\mathbf{y}}\otimes
\mathbf{M}_{\mathbf{x}}^{q}\rightarrow\mathbf{M}_{\left\vert \mathbf{x}%
\right\vert }^{\left\vert \mathbf{y}\right\vert },$ extend $\gamma$ to an
\emph{upsilon product }$\Upsilon:\overline{\mathbf{M}}\otimes\overline
{\mathbf{M}}\rightarrow\overline{\mathbf{M}}$ via
\begin{equation}
\Upsilon\left(  A\otimes B\right)  _{i\ell}=\left\{
\begin{array}
[c]{ll}%
\gamma\left(  A_{i\ell}^{\prime}\otimes B_{i\ell}^{\prime}\right)  , &
\text{if}\ A\otimes B\ \text{is a}\ \text{BTP}\\
& \\
0, & \text{otherwise.}%
\end{array}
\right.  \label{upsilon}%
\end{equation}
Then $\Upsilon$ sends a BTP $A^{q\times s}\otimes B^{t\times p}\in
\mathbf{M}_{\mathbf{v}}^{\mathbf{y}}\otimes\mathbf{M}_{\mathbf{x}}%
^{\mathbf{u}}$ with $A_{i\ell}^{\prime}\otimes B_{i\ell}^{\prime}\in
\mathbf{M}_{p_{\ell}}^{\mathbf{y}_{i}}\otimes\mathbf{M}_{\mathbf{x}_{\ell}%
}^{q_{i}}$ to a $t\times s $ monomial in $\mathbf{M}_{\left\vert
\mathbf{x}_{1}\right\vert ,\ldots,\left\vert \mathbf{x}_{s}\right\vert
}^{\left\vert \mathbf{y}_{1}\right\vert ,\ldots,\left\vert \mathbf{y}%
_{t}\right\vert }.$ We denote $A\cdot B=\Upsilon(A\otimes B);$ when $\left[
\theta_{j}\right]  \otimes\left[  \eta_{i}\right]  $ is a TP we denote
$\gamma(\theta_{1},\ldots,\theta_{q};\eta_{1},\ldots,\eta_{p})=\left(
\theta_{1}\otimes\cdots\otimes\theta_{q}\right)  \cdot\left(  \eta_{1}%
\otimes\cdots\otimes\eta_{p}\right)  $. As an arrow, $A\cdot B$ runs from the
initial point of $B$ to the terminal point of $A.$ Note that $\mathbf{M\cdot
M\subseteq M}$ so that $\Upsilon$ restricts to an upsilon product on
$\mathbf{M}.$

\begin{example}
\label{upsilonex}Continuing Example \ref{example1}, the action of $\Upsilon$
on $A^{4\times2}\otimes B^{2\times3}\in\mathbf{M}_{2,1}^{1,5,4,3}%
\otimes\mathbf{M}_{1,2,3}^{3,1}$ produces a $2\times2$ monomial in
$\mathbf{M}_{3,3}^{10,3}:$\vspace{0.2in}\newline\hspace*{0.1in}%
\setlength{\unitlength}{1.1mm}\linethickness{0.4pt}\begin{picture}(118.66,29.34)
\put(5.20,25.00){\makebox(0,0)[cc]{$\theta_{1,2}$}}
\put(5.20,18.00){\makebox(0,0)[cc]{$\theta_{5,2}$}}
\put(5.20,11.67){\makebox(0,0)[cc]{$\theta_{4,2}$}}
\put(5.20,4.67){\makebox(0,0)[cc]{$\theta_{3,2}$}}
\put(1.66,8.67){\dashbox{0.67}(6.33,20.00)[cc]{}}
\put(1.66,2.34){\dashbox{0.67}(6.33,4.67)[cc]{}}
\put(13.85,25.00){\makebox(0,0)[cc]{$\theta_{1,1}$}}
\put(13.85,18.00){\makebox(0,0)[cc]{$\theta_{5,1}$}}
\put(13.85,11.67){\makebox(0,0)[cc]{$\theta_{4,1}$}}
\put(13.85,4.67){\makebox(0,0)[cc]{$\theta_{3,1}$}}
\put(10.33,8.67){\dashbox{0.67}(6.33,20.00)[cc]{}}
\put(10.33,2.34){\dashbox{0.67}(6.33,4.67)[cc]{}}
\put(1.67,15.34){\oval(3.33,28.00)[l]}
\put(16.80,15.00){\oval(3.33,28.00)[r]} \
\put(20.83,15.77){\makebox(0,0)[cc]{$\cdot$}} \
\put(28.00,20.01){\makebox(0,0)[cc]{$\eta_{3,1}$}}
\put(34.58,20.01){\makebox(0,0)[cc]{$\eta_{3,2}$}}
\put(42.53,20.01){\makebox(0,0)[cc]{$\eta_{3,3}$}}
\put(24.67,17.67){\dashbox{0.67}(12.67,5.67)[cc]{}}
\put(39.33,17.67){\dashbox{0.67}(6.67,5.67)[cc]{}}
\put(28.00,12.34){\makebox(0,0)[cc]{$\eta_{1,1}$}}
\put(34.58,12.34){\makebox(0,0)[cc]{$\eta_{1,2}$}}
\put(42.53,12.34){\makebox(0,0)[cc]{$\eta_{1,3}$}}
\put(24.67,10.01){\dashbox{0.67}(12.67,5.67)[cc]{}}
\put(39.53,10.01){\dashbox{0.67}(6.67,5.67)[cc]{}}
\put(23.90,16.60){\oval(2.00,15.33)[l]}
\put(46.56,16.60){\oval(2.00,15.33)[r]}
\put(51.00,16.00){\makebox(0,0)[cc]{$=$}} \ \
\put(69.4,20.67){\makebox(0,0)[cc]{$_{\gamma(\theta_{1,2},\theta_{5,2},\theta_{4
,2};\,\eta_{3,1},\eta_{3,2})}$}}
\put(97.00,20.67){\makebox(0,0)[cc]{$_{\gamma(\theta_{1,1},\theta_{5,1},\theta_{
4,1};\,\eta_{3,3})}$}}
\put(69.31,12.67){\makebox(0,0)[cc]{$_{\gamma(\theta_{3,2};\,\eta_{1,1},\eta_{1,
2})}$}}
\put(97.00,12.67){\makebox(0,0)[cc]{$_{\gamma(\theta_{3,1};\,\eta_{1,3})}$}}
\put(55.16,16.60){\oval(1.67,15.33)[l]}
\put(108.33,16.60){\oval(2.67,15.33)[r]}
\put(112.00,16.67){\makebox(0,0)[cc]{$.$}}
\end{picture}\vspace{0.1in}\newline In the target, $\left(  \left\vert
\mathbf{x}_{1}\right\vert ,\left\vert \mathbf{x}_{2}\right\vert \right)
=\left(  1+2,3\right)  $ since $\left(  p_{1},p_{2}\right)  =\left(
2,1\right)  ;$ and $\left(  \left\vert \mathbf{y}_{1}\right\vert ,\left\vert
\mathbf{y}_{2}\right\vert \right)  =\left(  1+5+4,3\right)  $ since $\left(
q_{1},q_{2}\right)  =\left(  3,1\right)  .$ As an arrow, $A\cdot B$
initializes at $\left(  6,2\right)  $ and terminates at $\left(  2,13\right)
.$\vspace{0.1in}
\end{example}

The applications below relate to the following special case: Let $H$ be a
graded module over a commutative ring with unity and view $M=End(TH)$ as a
bigraded module via $M_{n,m}={Hom}\left(  H^{\otimes m},H^{\otimes n}\right)
.$ Then a $q\times p$ monomial $A\in\mathbf{M}_{\mathbf{x}}^{\mathbf{y}}$
admits a representation as an operator on $\mathbb{N}^{2}$ via%
\[
\left(  H^{\otimes\left\vert \mathbf{x}\right\vert }\right)  ^{\otimes
q}\approx\left(  H^{\otimes x_{1}}\otimes\cdots\otimes H^{\otimes x_{p}%
}\right)  ^{\otimes q}\overset{A}{\rightarrow}\left(  H^{\otimes y_{1}%
}\right)  ^{\otimes p}\otimes\cdots\otimes\left(  H^{\otimes y_{q}}\right)
^{\otimes p}%
\]%
\[
\overset{\sigma_{y_{1},p}\otimes\cdots\otimes\sigma_{y_{q},p}}{\longrightarrow
}\left(  H^{\otimes p}\right)  ^{\otimes y_{1}}\otimes\cdots\otimes\left(
H^{\otimes p}\right)  ^{\otimes y_{q}}\approx\left(  H^{\otimes p}\right)
^{\otimes\left\vert \mathbf{y}\right\vert },
\]
where $\left(  s,t\right)  \in\mathbb{N}^{2}$ is identified with $\left(
H^{\otimes s}\right)  ^{\otimes t}$ and $\sigma_{s,t}:\left(  H^{\otimes
s}\right)  ^{\otimes t}\overset{\approx}{\rightarrow}\left(  H^{\otimes
t}\right)  ^{\otimes s}$ is the canonical permutation of tensor factors
$\sigma_{q,p}:\left(  \left(  a_{11}\cdots a_{q1}\right)  \cdots\left(
a_{1p}\cdots a_{qp}\right)  \right)  $\ $\mapsto\left(  \left(  a_{11}\cdots
a_{1p}\right)  \cdots\left(  a_{q1}\cdots a_{qp}\right)  \right)  $. The
canonical structure map is
\begin{equation}
\gamma=\sum\gamma_{\mathbf{x}}^{\mathbf{y}}:\mathbf{M}_{p}^{\mathbf{y}}%
\otimes\mathbf{M}_{\mathbf{x}}^{q}\overset{\iota_{p}\otimes\iota_{q}%
}{\longrightarrow}\mathbf{M}_{pq}^{|\mathbf{y}|}\otimes\mathbf{M}%
_{|\mathbf{x}|}^{qp}\overset{\text{\textit{id}}\otimes\sigma_{q,p}^{\ast}%
}{\longrightarrow}\mathbf{M}_{pq}^{|\mathbf{y}|}\otimes\mathbf{M}%
_{|\mathbf{x}|}^{pq}\overset{\circ}{\longrightarrow}\mathbf{M}_{|\mathbf{x}%
|}^{|\mathbf{y}|}, \label{gamma}%
\end{equation}
where $\iota_{p}$ and $\iota_{q}$ are the canonical isomorphisms and
$\sigma_{q,p}^{\ast}$ is induced by $\sigma_{q,p}$ (c.f. \cite{Adams},
\cite{MSS}), \ induces a canonical \emph{associative} $\Upsilon$ product on
$\mathbf{M}$ whose action on matrices of double corollas typically produces a
matrix of non-planar graphs (see Figure 2).\vspace{0.4in}

\hspace*{0.8in}\setlength{\unitlength}{0.0015in}\begin{picture}
(152,-130)(202,-90) \thicklines
\put(-280,50){\line ( 0,-1){1200}}
\put(-280,50){\line ( 1,0){100}}
\put(-280,-1150){\line (1,0){100}}
\put(0,-150){\line ( 1,-1){100}}
\put(0,-150){\line ( -1,-1){100}}
\put(0,-150){\line( 0,1){150}}
\put(0,-150){\makebox(0,0){$\bullet$}}
\put(0,-550){\line ( 1,-1){100}}
\put(0,-550){\line ( -1,-1){100}}
\put(0,-550){\line( 0,1){150}}
\put(0,-550){\makebox(0,0){$\bullet$}}
\put(0,-950){\line ( 1,-1){100}}
\put(0,-950){\line ( -1,-1){100}}
\put(0,-950){\line( 0,1){150}}
\put(0,-950){\makebox(0,0){$\bullet$}}
\put(280,50){\line ( 0,-1){1200}}
\put(280,50){\line ( -1,0){100}}
\put(280,-1150){\line (-1,0){100}}
\put(500,-550){\makebox(0,0){$\bullet$}}
\put(700,-300){\line ( 0,-1){500}}
\put(700,-300){\line ( 1,0){100}}
\put(700,-800){\line (1,0){100}}
\put(900,-550){\line ( 1,1){100}}
\put(900,-550){\line ( -1,1){100}}
\put(900,-700){\line( 0,1){254}}
\put(900,-550){\makebox(0,0){$\bullet$}}
\put(1200,-550){\line ( 1,1){100}}
\put(1200,-550){\line ( -1,1){100}}
\put(1200,-700){\line( 0,1){254}}
\put(1200,-550){\makebox(0,0){$\bullet$}}
\put(1400,-300){\line ( 0,-1){500}}
\put(1400,-300){\line ( -1,0){100}}
\put(1400,-800){\line (-1,0){100}}
\put(1600,-550){\makebox(0,0){$=$}}
\put(1900,-250){\line ( 1,-1){345}}
\put(2300,-650){\line ( 1,-1){270}}
\put(1900,-250){\line ( -1,-1){100}}
\put(1800,-350){\line (0,-1){399}}
\put(1900,-250){\line( 0,1){150}}
\put(2270,-250){\line ( 1,-1){297}}
\put(2270,-250){\line ( -1,-1){165}}
\put(1970,-550){\line (1,1){93}}
\put(1970,-550){\line (0,-1){550}}
\put(2270,-250){\line( 0,1){150}}
\put(2640,-250){\line ( 1,-1){100}}
\put(2740,-350){\line (0,-1){399}}
\put(2640,-250){\line ( -1,-1){165}}
\put(2430,-460){\line ( -1,-1){460}}
\put(2640,-250){\line( 0,1){150}}
\put(1970,-920){\line(-1,1){170}}
\put(2568,-548){\line (0,-1){552}}
\put(2568,-920){\line(1,1){171}}
\put(2270,-250){\makebox(0,0){$\bullet$}}
\put(1900,-250){\makebox(0,0){$\bullet$}}
\put(2640,-250){\makebox(0,0){$\bullet$}}
\put(1970,-920){\makebox(0,0){$\bullet$}}
\put(2568,-920){\makebox(0,0){$\bullet$}}
\end{picture}\vspace*{1.7in}

\begin{center}
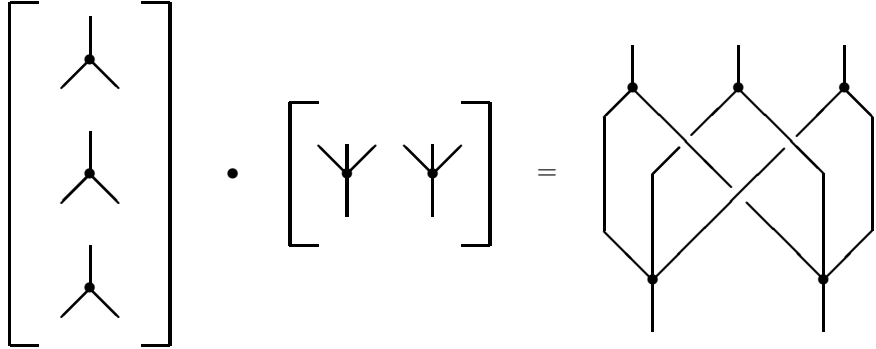
Figure 2. The $\gamma$-product as a non-planar graph.\vspace{0.2in}
\end{center}

\noindent In this setting, $\gamma$ agrees with the composition product on the
universal preCROC \cite{Shoikhet}.\vspace{0.1in}

\section{Cup products}

The two pairs of dual cup products defined in this section play an essential
role in the theory of structure relations. Let $\left(  H,d\right)  $ be a DG
module over a commutative ring with unity. For each $i,j\geq2,$ choose
operations $m_{i}:H^{\otimes i}\rightarrow H$ and $\Delta_{j}:H\rightarrow
H^{\otimes j}$ thought of as elements of $M=End\left(  TH\right)  .$ Recall
that planar rooted trees (PRT's) parametrize the faces of Stasheff's
associahedra $K=%
{\textstyle\bigsqcup\limits_{n\geq2}}
K_{n}\ $and provide module generators for cellular chains $C_{\ast}\left(
K\right)  $ \cite{MSS}. Whereas top dimensional faces correspond with
corollas, lower dimensional faces correspond with more general PRT's. Now
given a face $a\subseteq K,$ consider the class of all planar rooted trees
with levels (PLT's) representing $a$ and choose a representative with exactly
one node in each level. In this way, we obtain a particularly nice set of
module generators for $C_{\ast}\left(  K\right)  ,$ denoted by $\mathcal{K}$.
Note that the elements of a class of PLT's represent the same function
obtained by composing in various ways. The results obtained here are
independent of choice since they depend only on the function.

Let $G$ be a DGA concentrated in degree zero and consider the cellular
cochains on $K$ with coefficients in $G$:%
\[
C^{p}\left(  K;G\right)  =Hom^{-p}(C_{p}\left(  K\right)  ;G).
\]
A diagonal $\Delta$ on $C_{\ast}\left(  K\right)  $ induces a cup product
$\smile$ on $C^{\ast}\left(  K;G\right)  $ via%
\[
f\smile g=\cdot\left(  f\otimes g\right)  \Delta,
\]
where \textquotedblleft$\cdot$\textquotedblright\ denotes multiplication in
$G.$

The essential submodule $\mathbf{M,}$ which serves as our coefficient module,
is canonically endowed with dual associative \emph{wedge} and \emph{\v{C}ech
cross products }defined on a monomial pair $A\otimes B\in\mathbf{M}%
_{\mathbf{v}}^{\mathbf{y}}\otimes\mathbf{M}_{\mathbf{x}}^{\mathbf{u}}$ by
\[
A\overset{_{\wedge}}{\times}B=\left\{
\begin{array}
[c]{ll}%
A\otimes B, & \text{if }\mathbf{v}=\mathbf{x,}\\
0, & \text{otherwise,}%
\end{array}
\right.  \text{ \ and \ }A\overset{\vee}{\times}B=\left\{
\begin{array}
[c]{ll}%
A\otimes B, & \text{if }\mathbf{u}=\mathbf{y,}\\
0, & \text{otherwise.}%
\end{array}
\right.
\]
Denote $\overset{\wedge}{\mathbf{M}}=\left(  \mathbf{M,}\overset{_{\wedge}%
}{\times}\right)  $ and $\overset{\vee}{\mathbf{M}}=\left(  \mathbf{M,}%
\overset{_{\vee}}{\times}\right)  $ and note that $\mathbf{M}_{\mathbf{x}%
}^{\mathbf{y}}\overset{_{\wedge}}{\times}\mathbf{M}_{\mathbf{x}}^{\mathbf{u}%
}\subseteq\mathbf{M}_{\mathbf{x}}^{\mathbf{y,u}}$ and $\mathbf{M}_{\mathbf{v}%
}^{\mathbf{y}}\overset{_{\vee}}{\times}\mathbf{M}_{\mathbf{x}}^{\mathbf{y}%
}\subseteq\mathbf{M}_{\mathbf{v,x}}^{\mathbf{y}}.$ Thus non-zero cross
products concatenate matrices:%
\[
A\overset{_{\wedge}}{\times}B=%
\genfrac{[}{]}{0pt}{1}{A}{B}%
\text{ \ and \ }A\overset{_{\vee}}{\times}B=\left[  A\text{ }B\right]  .
\]
As arrows, $A\overset{_{\wedge}}{\times}B$ runs from vertical $x=\left\vert
\mathbf{x}\right\vert $ to vertical $x=p,$ whereas $A\overset{_{\vee}}{\times
}B$ runs from horizontal $y=q$ to horizontal $y=\left\vert \mathbf{y}%
\right\vert .$ In particular, if $A\in\mathbf{M}_{a}^{b}$, then $A^{\overset
{_{\wedge}}{\times}n}\in\mathbf{M}_{a}^{b\cdots b}$ is an arrow from $\left(
a,n\right)  $ to $\left(  1,nb\right)  $ and $A^{\overset{_{\vee}}{\times}%
n}\in\mathbf{M}_{a\cdots a}^{b}$ is an arrow from $\left(  na,1\right)  $ to
$\left(  n,b\right)  .$ These cross products together with the S-U diagonal
$\Delta_{K}$ \cite{SU2} induce wedge and \v{C}ech cup products $\wedge$ and
$\vee$ in $C^{\ast}(K;\overset{_{\wedge}}{\mathbf{M}})$ and $C^{\ast
}(K;\overset{_{\vee}}{\mathbf{M}}),$ respectively.

The modules $C^{\ast}(K;\overset{_{\wedge}}{\mathbf{M}})$ and $C^{\ast
}(K;\overset{_{\vee}}{\mathbf{M}})$ are equipped with second cup products
$\wedge_{\ell}$ and $\vee_{\ell}$ arising from the $\Upsilon$-product on
$\mathbf{M}$ together with the \textquotedblleft leaf
coproduct\textquotedblright\ $\Delta_{\ell}:C_{\ast}\left(  K\right)
\rightarrow C_{\ast}\left(  K\right)  \otimes C_{\ast}\left(  K\right)  ,$
which we now define. Let $T=T^{1}\in\mathcal{K}$ be a $k$-level PLT. Prune $T$
immediately below the first (top) level, trimming off a single corolla with
$n_{1}$ leaves and $r_{1}-1$ stalks. Numbering from left-to-right, let $i_{1}$
be the position of the corolla. The ($\emph{first}$) \emph{leaf sequence} of
$T$ is the $r_{1}$-tuple $\mathbf{x}_{i_{1}}\left(  n_{1}\right)  =\left(
1\cdots n_{1}\cdots1\right)  $ with $n_{1}$\ in position $i_{1}$ and $1$'s
elsewhere. Label the pruned tree $T^{2};$ inductively, the $\emph{j}^{th}%
$\emph{\ leaf sequence} of $T$ is the leaf sequence of $T^{j}.$ The induction
terminates when $j=k,$ in which case $i_{k}=r_{k}=1$ and $\mathbf{x}_{i_{k}%
}\left(  n_{k}\right)  =n_{k}.$ The \emph{descent sequence of }$T$ is the
$k$-tuple $\left(  \mathbf{x}_{i_{1}}\left(  n_{1}\right)  ,...,\mathbf{x}%
_{i_{k}}\left(  n_{k}\right)  \right)  .$

\begin{definition}
Let $T\in\mathcal{K}$ and identify $T$ with its descent sequence
$\mathbf{n=}\left(  \mathbf{n}_{1},...,\mathbf{n}_{k}\right)  .$ \textit{The
\underline{\textit{leaf coproduct}} of }$T$ is given by%
\[
\Delta_{\ell}\left(  T\right)  =\left\{
\begin{array}
[c]{cc}%
{\displaystyle\sum\limits_{2\leq i\leq k}}
\left(  \mathbf{n}_{1},...,\left\vert \mathbf{n}_{i}\right\vert \right)
\otimes\left(  \mathbf{n}_{i},\mathbf{n}_{i+1},...,\mathbf{n}_{k}\right)  , &
k>1\\
0, & k=1.
\end{array}
\right.
\]

\end{definition}

\noindent Define the \emph{leaf cup products} $\wedge_{\ell}$ and $\vee_{\ell
}$ on $C^{\ast}(K;\overset{\wedge}{\mathbf{M})}$ and $C^{\ast}(K;\overset
{\vee}{\mathbf{M})}$ by%
\[
f\wedge_{\ell}g=\cdot\left(  f\otimes g\right)  \tau\Delta_{\ell}\text{ and
}f\vee_{\ell}g=\cdot\left(  f\otimes g\right)  \Delta_{\ell},
\]
where $\tau$ interchanges tensor factors and $\cdot$ denotes the $\Upsilon$-product.

Note that all cup products defined in this section are non-associative and
non-commutative. Unless explicitly indicated otherwise, iterated cup products
are parenthesized on the extreme left, e.g., $f\vee g\vee h=\left(  f\vee
g\right)  \vee h.$

\section{Special $A_{\infty}$-bialgebras}

Structural compatibility of $d,$ the $m_{i}$'s and $\Delta_{j}$'s is expressed
in terms of the (restricted) biderivative $d_{\omega}$ and the
\textquotedblleft fraction product\textquotedblright\ $\bullet$ by the
equation $d_{\omega}\bullet d_{\omega}=0.$ We begin with a construction of the
biderivative in our restricted setting. Let $\varphi\in C^{\ast}%
(K;\overset{\wedge}{\mathbf{M}})$ and $\psi\in C^{\ast}(K;\overset{\vee
}{\mathbf{M})}$ be the cochains with top dimensional support such that
\[
\varphi\left(  e^{i-2}\right)  =m_{i}\text{ and }\psi\left(  e^{j-2}\right)
=\Delta_{j}.
\]
We think of $\varphi$ and $\psi$ as acting on uprooted and downrooted trees,
respectively (see Figure 3).\vspace{0.8in}\newline\hspace*{1.3in}%
\setlength{\unitlength}{0.0015in}\begin{picture}
(152,-130)(202,-90) \thicklines
\put(300,150){\line( -1,-1){150}}
\put(300,150){\line( 1,-1){150}}
\put(300,275){\line( 0,-1){280}}
\put(300,150){\line(1,-2){75}}
\put(300,150){\line( -1,-2){75}}
\put(300,142){\makebox(0,0){$\bullet$}}
\put(0,142){\makebox(0,0){$\varphi$}}
\put(100,300){\line( 0,-1){330}}
\put(100,300){\line( 1,0){30}}
\put(100,-30){\line( 1,0){30}}
\put(500,300){\line( 0,-1){330}}
\put(500,300){\line( -1,0){30}}
\put(500,-30){\line( -1,0){30}}
\put(680,142){\makebox(0,0){$=  m_{5}$ }}
\put(1600,125){\line( -1,1){150}}
\put(1600,125){\line( 1,1){150}}
\put(1600,275){\line( 0,-1){280}}
\put(1600,125){\line(1,2){75}}
\put(1600,125){\line( -1,2){75}}
\put(1600,130){\makebox(0,0){$\bullet$}}
\put(1300,142){\makebox(0,0){$\psi$}}
\put(1400,300){\line( 0,-1){330}}
\put(1400,300){\line( 1,0){30}}
\put(1400,-30){\line( 1,0){30}}
\put(1800,300){\line( 0,-1){330}}
\put(1800,300){\line( -1,0){30}}
\put(1800,-30){\line( -1,0){30}}
\put(1980,142){\makebox(0,0){$=  \Delta_5$}}
\end{picture}\vspace*{0.1in}

\begin{center}
Figure 3: The actions of $\varphi$ and $\psi.\vspace*{0.2in}$
\end{center}

Let $T^{c}H$ denote the tensor coalgebra of $H.$ The \emph{coderivation
cochain of }$\varphi$ is the cochain $\varphi^{c}\in C^{\ast}(K;\overset
{\wedge}{\mathbf{M}})$ that extends $\varphi$ to cells of $K$ in \textit{codim
}$1$ such that
\[%
{\displaystyle\sum\limits_{\text{\textit{codim} }e\text{ }=\text{ }0,1}}
\varphi^{c}\left(  e\right)  \in Coder\left(  T^{c}H\right)
\]
is the cofree linear coextension of $\varphi\left(  K\right)  =\sum_{i\geq
2}\varphi\left(  e^{i-2}\right)  $ as a coderivation. Thus if $T\in
\mathcal{K}$ is an uprooted $2$-level tree with $n+k$ leaves and leaf sequence
$\mathbf{x}_{i}\left(  k\right)  $,
\[
\varphi^{c}\left(  T\right)  =1^{\otimes i-1}\otimes m_{k}\otimes1^{\otimes
n-i+1}=\left[  1\text{ }\cdots\text{ }m_{k}\text{ }\cdots\text{ }1\right]
\in\mathbf{M}_{\mathbf{x}_{i}\left(  k\right)  }^{1}%
\]
and is represented by the arrow from $\left(  n+k,1\right)  $ to $\left(
n+1,1\right)  $ on the horizontal axis in $\mathbb{N}^{2}.$ Dually, let
$T^{a}\left(  H\right)  $ denote the tensor algebra of $H.$ The
\emph{derivation cochain of} $\psi$ is the cochain $\psi^{a}\in C^{\ast
}(K;\overset{\vee}{\mathbf{M})}$ that extends $\psi$ to cells of $K$ in
\textit{codim }$1$ such that
\[%
{\displaystyle\sum\limits_{\text{\textit{codim} }e\text{ }=\text{ }0,1}}
\psi^{a}\left(  e\right)  \in Der\left(  T^{a}H\right)
\]
is the free linear extension of $\psi\left(  K\right)  =\sum_{i\geq2}%
\psi\left(  e^{i-2}\right)  $ as a derivation. Thus if $T\in\mathcal{K}$ is an
downrooted $2$-level tree with $n+k$ leaves and leaf sequence $\mathbf{y}%
_{i}\left(  k\right)  $,
\[
\psi^{a}\left(  T\right)  =1^{\otimes i-1}\otimes\Delta_{k}\otimes1^{\otimes
n-i+1}=\left[  1\text{ }\cdots\text{ }\Delta_{k}\text{ }\cdots\text{
}1\right]  ^{T}\in\mathbf{M}_{1}^{\mathbf{y}_{i}\left(  k\right)  }%
\]
and is represented by the arrow from $\left(  1,n+1\right)  $ to $\left(
1,n+k\right)  $ on the vertical axis.

Evaluating leaf cup powers of $\varphi^{c}$ (respt. $\psi^{a}$) generates a
representative of each class of compositions involving the $m_{i}$'s (respt.
$\Delta_{j}$'s). So let
\[
\xi=\varphi^{c}+\varphi^{c}\wedge_{\ell}\varphi^{c}+\cdots+\left(  \varphi
^{c}\right)  ^{\wedge_{\ell}k}+\cdots\text{ }%
\]%
\[
\zeta=\psi^{a}+\psi^{a}\vee_{\ell}\psi^{a}+\cdots+\left(  \psi^{a}\right)
^{\vee_{\ell}k}+\cdots
\]
and note that if $e$ is a cell of $K,$ each non-zero component of $\xi\left(
e\right)  $ (respt. $\zeta\left(  e\right)  $) is represented by a
left-oriented horizontal (respt. upward-oriented vertical) arrow.

Furthermore, evaluating wedge and \v{C}ech cup powers of $\xi$ (respt. $\zeta
$) generates the components of the cofree coextension of $\xi\left(  K\right)
$ as a $\Delta_{K}$-coderivation (respt. free extension of $\zeta\left(
K\right)  $ as a $\Delta_{K}$-derivation). So let%
\[
\overset{\wedge}{\varphi}=\xi+\xi\wedge\xi+\cdots+\xi^{\wedge k}+\cdots
\]%
\[
\overset{\vee}{\psi}=\zeta+\zeta\vee\zeta+\cdots+\zeta^{\vee k}+\cdots
\]
and note that the component $\xi^{\wedge k}\left(  e^{i-2}\right)  :\left(
H^{\otimes i}\right)  ^{\otimes k}\rightarrow\left(  H^{\otimes1}\right)
^{\otimes k}$ is represented by a left-oriented horizontal arrow from $\left(
i,k\right)  $ to $\left(  1,k\right)  $ while the component $\zeta^{\vee
k}\left(  e^{i-2}\right)  :\left(  H^{\otimes1}\right)  ^{\otimes
k}\rightarrow\left(  H^{\otimes i}\right)  ^{\otimes k}$ is represented by a
upward-oriented vertical arrow from $\left(  k,1\right)  $ to $\left(
k,i\right)  .$

Let $M_{0}=M_{1,1}$. For reasons soon to become clear, the only structure
relations involving the differential $d$ are the classical quadratic relations
in an $A_{\infty}$-(co)algebra. Note that $d\in M_{0}$ and let $\mathbf{1}%
^{s}=\left(  1,\ldots,1\right)  \in\mathbb{N}^{s}.$ Given $\theta\in M_{0}$
and $p,q\geq1,$ consider the monomials $\theta_{i}^{q\times1}\in
\mathbf{M}_{\mathbf{1}}^{\mathbf{1}^{q}}$ and $\theta_{j}^{1\times p}%
\in\mathbf{M}_{\mathbf{1}^{p}}^{1}$ all of whose entries are the identity
except the $i^{th}$ in $\theta_{i}^{q\times1}$ and the $j^{th}$ in $\theta
_{j}^{1\times p},$ both of which are $\theta.$ Define $Bd_{0}:M_{0}%
\rightarrow\mathbf{M}$ by%
\[
Bd_{0}(\theta)=\sum\limits_{\substack{1\leq i\leq q,\text{ }1\leq j\leq p
\\p,q\geq1}}\theta_{i}^{q\times1}+\theta_{j}^{1\times p}.
\]
Then $Bd_{0}\left(  \theta\right)  $ is the (co)free linear (co)extension of
$\theta$ as a (co)derivation. Note that each component of $Bd_{0}\left(
\theta\right)  $ is represented by an arrow of \textquotedblleft
length\textquotedblright\ zero.

Let $M_{1}=\left(  M{_{1,\ast}}\oplus M_{\ast,1}\right)  /M_{1,1}$\ and define
$Bd_{1}:M_{1}\rightarrow\mathbf{M}$ by%
\begin{equation}
Bd_{1}\left(  \theta\right)  =%
{\displaystyle\sum\limits_{\substack{e\text{ }\subseteq\text{ }K \\\text{codim
}e\text{ }=\text{ }0}}}
(\overset{\wedge}{\varphi}+\overset{\vee}{\psi})\left(  e\right)  +%
{\displaystyle\sum\limits_{\substack{e\text{ }\subseteq\text{ }K
\\\text{\textit{codim} }e\text{ }=\text{ }1}}}
\left(  \varphi^{c}+\psi^{a}\right)  \left(  e\right)  . \label{bider}%
\end{equation}
Note that the components of $Bd_{1}\left(  \theta\right)  $ are represented by
upward-oriented vertical arrows and left-oriented horizontal arrows; the
right-hand component of (\ref{bider}) is given by Gerstenhaber's $\circ$-(co)operation.

Let $\rho_{0}:\mathbf{M\rightarrow M}_{0}$ and $\rho_{1}:\mathbf{M\rightarrow
M}_{1}$ denote the canonical projections.

\begin{definition}
The \underline{restricted biderivative} is the (non-linear) map $d_{\underline
{\ }}:\mathbf{M}\rightarrow\mathbf{M}$ given by%
\[
d_{\underline{\ }}=Bd_{0}\circ\rho_{0}+Bd_{1}\circ\rho_{1}.
\]
The symbol $d_{\theta}$ denotes the restricted biderivative of $\theta.$
\end{definition}

\noindent Finally, the composition
\[
\bullet:\mathbf{M}\times\mathbf{M}\overset{d_{\underline{\ }}\otimes
d_{\underline{\ }}}{\longrightarrow}\mathbf{M}\times\mathbf{M}\overset
{\Upsilon}{\longrightarrow}\mathbf{M}%
\]
defines the \emph{fraction product}. Special $A_{\infty}$-bialgebras are
defined in terms of the fraction product as follows:

\begin{definition}
Let $\omega=d+%
{\textstyle\sum\nolimits_{i,j\geq2}}
\left(  m_{i}+\Delta_{j}\right)  \in M_{0}\oplus M_{1}.$ Then $\left(
H,d,m_{i},\Delta_{j}\right)  _{i,j\geq2}$ is a \underline{special $A_{\infty}%
$-bialgebra} provided
\[
d_{\omega}\bullet d_{\omega}=0.
\]

\end{definition}

\noindent Note that one recovers the classical quadratic relations in an
$A_{\infty}$-algebra when $\omega=d+%
{\textstyle\sum\nolimits_{i\geq2}}
m_{i}.$

\section{Structure Relations}

The structure relations in a special $A_{\infty}$-bialgebra $\left(
H,d,m_{i},\Delta_{j}\right)  _{i,j\geq2}$ follow easily from the following two observations:

\begin{enumerate}
\item If $\theta,\eta\in\mathbf{M,}$ then $\theta\bullet\eta=0$ whenever the
projection of $\theta$ or $\eta$ to $M_{0}\oplus M_{1}$ is zero.

\item Each non-zero component in the projections of $\theta$ and $\eta$ is
represented by a horizontal, vertical or zero length arrow.
\end{enumerate}

\noindent By (1), each component of $d_{\omega}\bullet d_{\omega}$ is a
\textquotedblleft transgression\textquotedblright\ represented by a
\textquotedblleft2-step\textquotedblright\ path of arrows from the horizontal
axis $M_{1,\ast}$ to the vertical axis $M_{\ast,1};$ and by (2), each such
2-step path follows the edges of a (possibly degenerate) rectangle positioned
with one of its vertices at $\left(  1,1\right)  $.

Now relations involving $d$ arise from degenerate rectangles since arrows of
length zero represent components in the (co)extensions of $d$. Hence $d$
interacts with the $m_{i}$'s or the $\Delta_{j}$'s exclusively and the
relations involving $d$ are exactly the classical quadratic relations in an
$A_{\infty}$-(co)algebra.

On the other hand, relations involving the $m_{i}$'s and $\Delta_{j}$'s arise
from non-degenerate rectangles since $m_{i}$ and $\Delta_{j}$ are represented
by the arrows $\left(  i,1\right)  \rightarrow\left(  1,1\right)  $ and
$\left(  1,1\right)  \rightarrow\left(  1,j\right)  $. While the two-step path
$\left(  i,1\right)  \rightarrow\left(  1,1\right)  \rightarrow\left(
1,j\right)  $ represents the (usual) composition $\Delta_{j}\bullet m_{i},$
the two-step path $\left(  i,1\right)  \rightarrow\left(  i,j\right)
\rightarrow\left(  1,j\right)  $ represents $\xi^{j}\left(  e^{i-2}\right)
\bullet\zeta^{i}\left(  e^{j-2}\right)  .$ Thus we obtain the relation%
\[
\Delta_{j}\bullet m_{i}=\xi^{j}\left(  e^{i-2}\right)  \bullet\zeta^{i}\left(
e^{j-2}\right)  .
\]
For example, by setting $i=j=2$ we obtain the classical bialgebra relation%
\[
\Delta_{2}\bullet m_{2}=\left[
\begin{array}
[c]{c}%
m_{2}\\
m_{2}%
\end{array}
\right]  \bullet\left[  \Delta_{2}\text{ }\Delta_{2}\right]  .
\]
And with $\left(  i,j\right)  =\left(  3,2\right)  $ we obtain
\[
\Delta_{2}\bullet m_{3}=\left\{  \left[
\begin{array}
[c]{c}%
m_{3}\\
m_{2}\left(  1\otimes m_{2}\right)
\end{array}
\right]  +\left[
\begin{array}
[c]{c}%
m_{2}\left(  m_{2}\otimes1\right) \\
m_{3}%
\end{array}
\right]  \right\}  \bullet\left[  \Delta_{2}\text{ }\Delta_{2}\text{ }%
\Delta_{2}\right]
\]
(see Figure 4).

\ \vspace*{0.5in}\newline\hspace*{1.2in}%
\setlength{\unitlength}{0.0015in}\begin{picture}
(152,-130)(202,-90) \thicklines
\put(2000,0){\vector( -1,0){2000}}
\put(2000,-1000){\vector( -1,0){2000}}
\put(2000,-1000){\vector( 0,1){1000}}
\put(1000,-1000){\vector( 0,1){1000}}
\put(0,-1000){\vector( 0,1){1000}}
\put(-200,-500){\makebox(0,0){$\left[ \Delta_{2} \right]$}}
\put(2400,-500){\makebox(0,0){$\left[ \Delta_{2} \hspace{.1in}  \Delta_{2}  \hspace{.1in}  \Delta_{2}  \right]$}}
\put(1000,-1200){\makebox(0,0){$\left[ m_{3} \right] $}}
\put(1000,200){\makebox(0,0){
$  \left[
\begin{array}{c}
m_{3} \\
m_{2}\left( 1\otimes m_{2}\right)
\end{array}%
\right] +\left[
\begin{array}{c}
m_{2}\left( m_{2}\otimes 1\right)  \\
m_{3}%
\end{array}%
\right]  $
}}
\put(1300,-500){\makebox(0,0){$\left[ \Delta_{2} \hspace{.1in}  \Delta_{2}    \right]$}}
\put(500,-850){\makebox(0,0){$\left[ m_{2} \right]$}}
\put(500,-175){\makebox(0,0){ $\left[
\begin{array}{c}
m_{2} \\
m_{2}%
\end{array}%
\right] $
}}
\end{picture}\vspace{1.9in}

\begin{center}
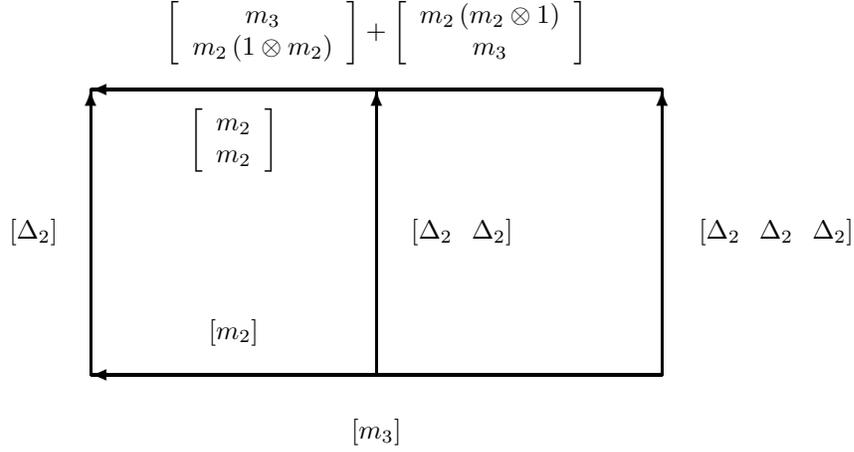
Figure 4: Some low order arrows in $\mathbf{M}.$\vspace{0.2in}
\end{center}

We summarize the discussion above in our main theorem:

\begin{theorem}
$\left(  H,d,m_{i},\Delta_{j}\right)  _{i,j\geq2}$ is a special $A_{\infty}%
$-bialgebra if $\left(  H,d,m_{i}\right)  _{i\geq2}$ is an $A_{\infty}%
$-algebra, $\left(  H,d,\Delta_{j}\right)  _{j\geq2}$ is an $A_{\infty}%
$-coalgebra and for all $i,j\geq2,$
\[
\Delta_{j}\bullet m_{i}=\xi^{j}\left(  e^{i-2}\right)  \bullet\zeta^{i}\left(
e^{j-2}\right)  .
\]
\vspace*{0.2in}\ 
\end{theorem}

\end{document}